\documentclass[]{aspm}

\articleinfo{}{}{}

\usepackage{verbatim}
\usepackage{amssymb}
\usepackage{amsbsy}
\usepackage{amscd}
\usepackage{amsmath}
\usepackage{amsthm}
\usepackage[mathscr]{eucal}

\newtheorem{thm}{Theorem}
\numberwithin{defn}{section}

\numberwithin{equation}{section}

\newtheorem{lem}{Lemma}[section]

\newtheorem{rem}{Remark}[section]

\def\R{{\bf R}}

\def\N{{\bf N}}

\def\d{\displaystyle}
\def\e{{\varepsilon}}

\def\wt{\widetilde}

\def\p{\partial}

\title[Semilinear Wave Equations in Two Space Dimensions]
{The sharp lower bound of the lifespan of solutions\\
to semilinear wave equations with low powers\\
in two space dimensions
\footnote{Dedicated to Professor Nakao Hayashi on the occasion of 
his sixties birthday.}}

\author[T.Imai]{Takuto Imai}
\author[M.Kato]{Masakazu Kato}
\author[H.Takamura]{Hiroyuki Takamura}
\author[K.Wakasa]{Kyouhei Wakasa}

\address{The first year of the master course,
Graduate School of Systems Information Science,
Future University Hakodate,
116-2 Kamedanakano-cho,
Hakodate, Hokkaido 041-8655, Japan.}
\address{College of Liberal Arts, Mathematical Science Research Unit, 
Muroran Institute of Technology, 27-1, Mizumoto-cho, 
Muroran, Hokkaido 050-8585, Japan.}
\address{Department of Complex and Intelligent Systems,
Faculty of Systems Information Science,
Future University Hakodate,
116-2 Kamedanakano-cho,
Hakodate, Hokkaido 041-8655, Japan.}
\address{College of Liberal Arts, Mathematical Science Research Unit, 
Muroran Institute of Technology, 27-1, Mizumoto-cho, 
Muroran, Hokkaido 050-8585, Japan.}
\address{}

\email{g2116004@fun.ac.jp}
\email{mkato@mmm.muroran-it.ac.jp}
\email{takamura@fun.ac.jp}
\email{wakasa@mmm.muroran-it.ac.jp}

\rcvdate{}
\rvsdate{}



\begin{document}

\begin{abstract}
This paper is devoted to a proof of the conjecture in Takamura \cite{Ta15}
on the lower bound of the lifespan of solutions to semilinear wave equations
in two space dimensions.
The result is divided into two cases according to the total integral of the initial speed. 
\end{abstract}

\maketitle

\section{Introduction}
\par
We consider the initial value problem,
\begin{equation}
\label{IVP}
\left\{
\begin{array}{ll}
u_{tt}-\Delta u=|u|^p
&\mbox{in}\quad \R^n\times[0,\infty),\\
u(x,0)=\e f(x),\ u_t(x,0)=\e g(x),
& x\in\R^n,
\end{array}
\right.
\end{equation}
where $u=u(x,t)$ is an unknown function,
$f$ and $g$ are given smooth functions of compact support
and $\e>0$ is \lq\lq small."
Let us define a lifespan $T(\e)$ of a solution of (\ref{IVP}) by
\[
T(\e):=\sup\{t>0:\exists\ \mbox{a solution $u$ of (\ref{IVP})
for arbitrarily fixed $(f,g)$}\},
\]
where \lq\lq solution" means a classical one when $p\ge2$.
When $1<p<2$, it means a weak one, but sometimes the one
given by associated integral equations to (\ref{IVP})
by standard Strichartz's estimate.
See Sideris \cite{Si84} for instance. 
\par
When $n=1$, we have $T(\e)<\infty$ for any power $p>1$ by Kato \cite{Kato80}.
When $n\ge2$, we have the following Strauss' conjecture on (\ref{IVP})
by Strauss \cite{St81}. 
\[
\begin{array}{lll}
T(\e)=\infty & \mbox{if $p>p_0(n)$ and $\e$ is \lq\lq small"}
& \mbox{(global-in-time existence)},\\
T(\e)<\infty & \mbox{if $1<p\le p_0(n)$}
& \mbox{(blow-up in finite time)},
\end{array}
\]
where $p_0(n)$ is so-called Strauss' exponent defined
by positive root of the quadratic equation, $\gamma(p,n)=0$, where
\begin{equation}
\label{gamma}
\gamma(p,n):=2+(n+1)p-(n-1)p^2.
\end{equation}
That is,
\begin{equation}
\label{p_0(n)}
p_0(n)=\frac{n+1+\sqrt{n^2+10n-7}}{2(n-1)}. 
\end{equation}
We note that $p_0(n)$ is monotonously decreasing in $n$.
This conjecture had been verified by many authors with partial results.
All the references on the final result in each part
can be summarized in the following table.
\begin{center}
\begin{tabular}{|c||c|c|c|}
\hline
& $p<p_0(n)$ & $p=p_0(n)$ & $p>p_0(n)$ \\
\hline
\hline
$n=2$ & Glassey \cite{G81a} & Schaeffer \cite{Sc85} & Glassey \cite{G81b}\\
\hline
$n=3$ & John \cite{J79} & Schaeffer \cite{Sc85} & John \cite{J79}\\
\hline
$n\ge4$ & Sideris \cite{Si84} &  
$
\begin{array}{l}
\mbox{Yordanov $\&$ Zhang \cite{YZ06}}\\
\mbox{Zhou \cite{Z07}, indep.}
\end{array}
$
&
$
\begin{array}{l}
\mbox{Georgiev $\&$ Lindblad}\\
\mbox{$\&$ Sogge \cite{GLS97}}\\
\end{array}
$
\\
\hline
\end{tabular} 
\end{center}
\par
In the blow-up case, i.e. $1<p\le p_0(n)$,
we are interested in the estimate of the lifespan $T(\e)$.
From now on, $c$ and $C$ stand for positive constants but independent of $\e$.
When $n=1$, we have the following estimate of the lifespan $T(\e)$ for any $p>1$.
\begin{equation}
\label{lifespan_1d}
\left\{
\begin{array}{cl}
c\e^{-(p-1)/2}\le T(\e)\le C\e^{-(p-1)/2}
& \mbox{if}\quad\d\int_{\R}g(x)dx\neq0,\\
c\e^{-p(p-1)/(p+1)}\le T(\e)\le C\e^{-p(p-1)/(p+1)}
&\mbox{if} \quad\d\int_{\R}g(x)dx=0.
\end{array}
\right.
\end{equation}
This result has been obtained by Zhou \cite{Z92_one}.
Moreover, Lindblad \cite{L90} has obtained more precise result for $p=2$,
\begin{equation}
\label{lifespan_1d_lim}
\left\{
\begin{array}{ll}
\d \exists \lim_{\e\rightarrow+0}\e^{1/2}T(\e)>0
&\mbox{if}\quad\d\int_{\R}g(x)dx\neq0,\\
\d \exists \lim_{\e\rightarrow+0}\e^{2/3}T(\e)>0
&\mbox{if}\quad\d\int_{\R}g(x)dx=0.
\end{array}
\right.
\end{equation}
Similarly to this, Lindblad \cite{L90} has also obtained the following result
for $(n,p)=(2,2)$.
\begin{equation}
\label{lifespan_2d_lim}
\left\{
\begin{array}{ll}
\d \exists \lim_{\e\rightarrow+0}a(\e)^{-1}T(\e)>0
&\mbox{if}\quad\d\int_{\R^2}g(x)dx\neq0\\
\d \exists \lim_{\e\rightarrow+0}\e T(\e)>0
&\mbox{if}\quad\d\int_{\R^2}g(x)dx=0,
\end{array}
\right.
\end{equation}
where $a=a(\e)$ is a number satisfying 
\begin{equation}
\label{a}
a^2\e^2\log(1+a)=1.
\end{equation}
When $1<p<p_0(n)\ (n\ge3)$ or $2<p<p_0(2)\ (n=2)$,
we have the following conjecture.
\begin{equation}
\label{lifespan_high-d}
c\e^{-2p(p-1)/\gamma(p,n)}\le T(\e)\le C\e^{-2p(p-1)/\gamma(p,n)},
\end{equation}
where $\gamma(p,n)$ is defined by (\ref{gamma}).
We note that (\ref{lifespan_high-d}) coincides
with the second line in (\ref{lifespan_1d})
if we define $\gamma(p,n)$ by (\ref{gamma}) even for $n=1$.
All the results verifying this conjecture are summarized in the following table.
\begin{center}
\begin{tabular}{|c||c|c|c|}
\hline
& lower bound of $T(\e)$ & upper bound of $T(\e)$\\
\hline
\hline
$n=2$ & Zhou \cite{Z93} & Zhou \cite{Z93}\\
\hline
$n=3$ & Lindblad \cite{L90} & Lindblad \cite{L90} \\
\hline  
$n\ge4$ 
&
Lai $\&$ Zhou \cite{LZ14}
&
Takamura \cite{Ta15}
\\
\hline
\end{tabular} 
\end{center}
We note that, for $n=2,3$,
\[
 \exists \lim_{\e\rightarrow+0}\e^{2p(p-1)/\gamma(p,n)}T(\e)>0
\]
is established in this table.
When $p=p_0(n)$,
we have the following conjecture.
\begin{equation}
\label{lifespan_critical}
\exp\left(c\e^{-p(p-1)}\right)\le T(\e)\le\exp\left(C\e^{-p(p-1)}\right).
\end{equation}
All the results verifying this conjecture are also summarized in the following table.
\begin{center}
\begin{tabular}{|c||c|c|c|}
\hline
& lower bound of $T(\e)$ & upper bound of $T(\e)$\\
\hline
\hline
$n=2$ & Zhou \cite{Z93} & Zhou \cite{Z93}\\
\hline
$n=3$ & Zhou \cite{Z92_three} & Zhou \cite{Z92_three} \\
\hline  
$n\ge4$ 
&
$
\begin{array}{l}
\mbox{Lindblad $\&$ Sogge \cite{LS96}}\\
\mbox{: $n\le 8$ or radially symm. sol.}
\end{array}
$
&
Takamura $\&$ Wakasa \cite{TW11}\\
\hline
\end{tabular} 
\end{center}
\par
In this paper, we are interested in the case of the open part,
$n=2$ and $1<p<2$.
There is no conjecture before Remark 4.1 in Takamura \cite{Ta15}
in which the following is implicitly described. 
\begin{equation}
\label{lifespan_2d_low}
\left\{
\begin{array}{ll}
\d\exists\lim_{\e\rightarrow+0}T(\e)\e^{(p-1)/(3-p)}>0
& \mbox{if}\quad\d\int_{\R^2}g(x)dx\neq0,\\
\d\exists\lim_{\e\rightarrow+0}T(\e)\e^{2p(p-1)/\gamma(p,2)}>0
&\mbox{if} \quad\d\int_{\R^2}g(x)d x=0.
\end{array}
\right.
\end{equation}
Theorem 3.2 and 4.1 in Takamura \cite{Ta15} 
are the partial result of (\ref{lifespan_2d_low}),
the upper bounds of $T(\e)$.
Our result is devoted to the lower bounds as follows.

\begin{thm}
\label{thm:main1}
Let $n=2$, $1<p<2$ and $(f,g)\in C_0^3(\R^2)\times C_0^2(\R^2)$.
Then, there exists a positive constant $\e_0
=\e_0(f,g,p,k)$ such that the lifespan $T(\e)$ 
of solutions of (\ref{IVP}) satisfies that
\begin{equation}
\label{lower_lifespan}
\left\{
\begin{array}{ll}
T(\e)\ge c\e^{-(p-1)/(3-p)}
& \mbox{if}\quad\d\int_{\R^2}g(x)dx\neq0,\\
T(\e)\ge c\e^{-2p(p-1)/\gamma(p,2)}
&\mbox{if} \quad\d\int_{\R^2}g(x)d x=0
\end{array}
\right.
\end{equation}
for $0<\e\le\e_0$, where $c$ is a positive constant independent of $\e$.
\end{thm}

\par
This paper is organized as follows.
In the next section, we employ the linear decay estimate
and basic lemmas for a priori estimates.
In the third section, we prove a priori estimates.
The proof of the Theorem \ref{thm:main1} is in the final section.


\section{Preliminaries}
\label{section:Pre}
\par
Throughout this paper, we may assume that
$(f,g)\in C_0^3(\R^2)\times C_0^2(\R^2)$ satisfy
\begin{equation}
\label{supp_data}
{\rm supp}\ (f,g) \subset\{x\in\R^2\ :\ |x|\le k\},\ k>1.
\end{equation}
Set
 \begin{equation}
\label{u_L}
\begin{array}{ll}
\d u_L(x,t):=\frac{\partial}{\partial t}R(f|x,t)+R(g|x,t),\\
\d R(\phi|x,t):=\frac{1}{2\pi}\int_{|x-y|\le t}\frac{\phi(y)}{\sqrt{t^2-|x-y|^2}}dy
=\frac{t}{2\pi}\int_{|\xi|\le1}\frac{\phi(x+t\xi)}{\sqrt{1-|\xi|^2}}d\xi.
\end{array}
\end{equation}
Then, we note that $u_L$ satisfies that
\[
\left\{
\begin{array}{ll}
\d (u_L)_{tt}-\Delta u_L=0
& \mbox{in}\ \R^2\times[0,\infty),\\
u_L(x,0)=f(x),\ (u_L)_t(x,0)=g(x),
& x\in \R^2
\end{array}
\right.
\]
in the classical sense, and also that
\begin{equation}
\label{supp_u_L}
{\rm supp}\ u_L\subset\{(x,t)\in\R^2\times[0,\infty)\ :\ |x|\le t+k\}.
\end{equation}
We shall employ the following key lemma.

\begin{lem}[Lindblad \cite{L90}]
\label{lem:decay_est_v}
Let $u_L$ be the one in (\ref{u_L}).
Then, there exist positive constants
$C_{0}=C_0(\|f\|_{W^{3,1}(\R^2)},\|g\|_{W^{2,1}(\R^2)},k)$ 
and $\wt{C_0}=\wt{C_0}(k)$ such that $u_L$ satisfies
\begin{equation}
\label{decay_est}
\begin{array}{ll}
\d \sum_{|\alpha|\le1}|\nabla_x^\alpha u_L(x,t)|
\le
&\d\frac{\d \wt{C_{0}}\left|\int_{\R^2}g(x)dx\right|}
{(t+|x|+2k)^{1/2}(t-|x|+2k)^{1/2}}\\
&\d +\frac{C_0}{(t+|x|+2k)^{1/2}(t-|x|+2k)^{3/2}}
\end{array}
\end{equation}
in $\R^2\times[0,\infty)$.
\end{lem}
\begin{rem}
This is not exactly Lemma 7.1 in \cite{L90},
but is basically its refined version.
For the sake of completeness, we prove it.
\end{rem}
\par\noindent
{\bf Proof of Lemma \ref{lem:decay_est_v}.}
First we note that it is sufficient to show (\ref{decay_est})
for $\alpha=0$ because $\nabla_x^\alpha$ passes to
the integrand in the representation (\ref{u_L}).
Moreover, due to von Wahl \cite{W71}, or Klainerman \cite{Kl80}
as described in Glassey \cite{G81b},
we have that
\[
|u_L(x,t)|\le\frac{C}{\sqrt{1+t}}
\left(\|f\|_{W^{2,1}(\R^2)}+\|g\|_{W^{1,1}(\R^2)}\right)
\quad\mbox{in}\ \R^2\times[0,\infty),
\]
where $C$ is a positive constant independent of $f$ and $g$.
Therefore (\ref{decay_est}) is obtained by (\ref{supp_u_L})
for $-k\le t-|x|\le 2k$, or $t\le4k$.
\par
From now on, we are concentrated in the case of $\alpha=0$, $t-|x|\ge 2k$ and $t\ge4k$.
Set $r:=|x|$.
First we prove (\ref{decay_est}) in the interior domain,
\[
D_{\rm int}:=\{(x,t)\in\R^2\times[0,\infty)\ :\ t\ge2r,\ t\ge4k\}.
\]
Since
\[
|x-y|\le r+|y|\le \frac{t}{2}+k\le t
\quad\mbox{for $(x,t)\in D_{\rm int}$ and $|y|\le k$},
\]
we can rewrite $R(g|x,t)$ in (\ref{u_L}) as
\[
R(g|x,t)=\frac{1}{2\pi}\int_{\R^2}\frac{g(y)}{\sqrt{t^2-|x-y|^2}}dy
\quad\mbox{in}\ D_{\rm int}.
\]
We shall consider the following estimate:
\begin{equation}
\label{diff_est}
\begin{array}{l}
\d \left|2\pi R(g|x,t)
-\frac{1}{\sqrt{(t+r)(t-r)}}\int_{\R^2}g(y)dy\right|\\
\d\le\frac{1}{\sqrt{(t+r)(t-r)}}
\int_{\R^2}\frac{\left|h(x,y,t)\right|}{\sqrt{t^2-|x-y|^2}}|g(y)|dy,
\end{array}
\end{equation}
where
\[
h(x,y,t):=\sqrt{(t+r)(t-r)}-\sqrt{t^2-|x-y|^2}.
\]
Making use of Taylor expansion in $y$ at the origin, we get
\[
h(x,y,t)=\frac{-<x,y>+\theta|y|^2}{\sqrt{t^2-|x-\theta y|^2}}
\quad \mbox{with}\ 0<\theta<1.
\]
It follows from
\[
|-<x,y>+\theta |y|^2|\le (t/2+k)k
\]
and
\[
\begin{array}{l}
t-|x-\theta y|\ge t-(r+|y|)\ge t-r-k\ \ge t/2-k\ge t/4,\\
t+|x-\theta y|\ge t,\\
\end{array}
\]
for $(x,t)\in D_{\rm int}$ and $|y|\le k$ that
\[
|h(x,y,t)|\le (1+2k/t)k\le2k\quad \mbox{for} \ (x,t)\in D_{\rm int}\ \mbox{and}\ |y|\le k.
\]
Thus, the right hand-side of (\ref{diff_est}) is dominated by
\[
\frac{2k}{\sqrt{(t+r)(t-r)}}
\int_{\R^2}\frac{|g(y)|dy}{\sqrt{t^2-|x-y|^2}}.
\]
Similarly to the above, it follows from
\[\begin{array}{l}
t-|x-y|\ge t-r-k\ge t-r/2-r/2-t/4\ge(t-r)/2,\\
t+|x-y|\ge t\ge t/2+r\ge(t+r)/2,
\end{array}
\]
for $(x,t)\in D_{\rm int}$ and $|y|\le k$ that
\[
\left|2\pi R(g|x,t)
-\frac{1}{\sqrt{(t+r)(t-r)}}\int_{\R^2}g(y)dy\right|
\le\frac{4k\|g\|_{L^1(\R^2)}}{(t+r)(t-r)}
\]
in $D_{\rm int}$. 
We also obtain that
\[
\left|\frac{\partial}{\partial t}R(f|x,t)\right|
\le\frac{1}{\pi}\int_{\R^2}\frac{t|f(y)|}{(t^2-|x-y|^2)^{3/2}}dy
\le\frac{4\|f\|_{L^1(\R^2)}}{\pi\sqrt{t+r}(t-r)^{3/2}}
\]
in $D_{\rm int}$. 
Since $3(t-r)\ge t+r$ holds in $D_{\rm int}$, 
summing up all the estimates, we have the desired estimate
\[
|u_L(x,t)|\le
\frac{\d\left|\int_{\R^2}g(y)dy\right|}{2\pi\sqrt{(t+r)(t-r)}}
+\frac{4\sqrt{3}\|f\|_{L^1(\R^2)}+2k\|g\|_{L^1(\R^2)}}{\pi(t+r)(t-r)}
\]
in $D_{\rm int}$.
\par
Next we prove (\ref{decay_est}) in the exterior domain,
\[
D_{\rm ext}:=\{(x,t)\in \R^2\times[0,\infty)\ :\ r+2k\le t \le 2r\}.
\]
Here we employ the different representation formula from (\ref{u_L}),
\[
R(g|x,t)
=\frac{1}{2\sqrt{2}\pi\sqrt{r}}
\int_{\rho-\rho^2z/2}^\infty\frac{ds}{\sqrt{s-\rho+\rho^2z/2}}
\int_{s=<\omega,y>-|y|^2z/2}g(y)dS_y,
\]
where $\omega:=x/r\in S^1$, $\rho:=r-t$ and $z:=1/r$.
This is established by (6.2.4) in H\"ormander \cite{H97}.
Due to $z\le 1/(2k)$ in $D_{\rm ext}$ and $|y|\le k$ by (\ref{supp_data}), we get
\[
|s|\le|y|+\frac{|y|^2z}{2}
\le k+\frac{k^2}{2}\cdot\frac{1}{2k}=\frac{5}{4}k.
\]
Since
\[
\begin{array}{l}
\d2\sqrt{2}\pi\sqrt{r}R(g|x,t)-\frac{1}{\sqrt{-\rho+\rho^2z/2}}\int_{\R^2}g(x)dx\\
\d=\int_{-5k/4}^{5k/4}\frac{ds}{\sqrt{s-\rho+\rho^2z/2}}
\int_{s=<\omega,y>-|y|^2z/2}g(y)dS_y\\
\d\quad-\frac{1}{\sqrt{-\rho+\rho^2z/2}}\int_{-5k/4}^{5k/4}ds
\int_{s=<\omega,y>-|y|^2z/2}g(y)dS_y
\end{array}
\]
and
\[
\begin{array}{l}
\d\left|\frac{1}{\sqrt{s-\rho+\rho^2z/2}}
-\frac{1}{\sqrt{-\rho+\rho^2z/2}}\right|
=\left|\frac{\sqrt{-\rho+\rho^2z/2}-\sqrt{s-\rho+\rho^2z/2}}
{\sqrt{(s-\rho+\rho^2z/2)(-\rho+\rho^2z/2)}}\right|\\
\d\le\frac{|s|}
{\sqrt{(s+|\rho|)|\rho|}(\sqrt{-\rho+\rho^2z/2}+\sqrt{s-\rho+\rho^2z/2})}
\le\frac{|s|}{|\rho|\sqrt{(s+|\rho|)}}
\end{array}
\]
hold, we have that
\[
\begin{array}{l}
\d\left|2\sqrt{2}\pi\sqrt{r}R(g|x,t)
-\frac{1}{\sqrt{-\rho+\rho^2z/2}}\int_{\R^2}g(x)dx\right|\\
\d\le\frac{5k}{4|\rho|\sqrt{(|\rho|-5k/4)}}
\|g\|_{L^1(\R^2)}.
\end{array}
\]
We note that $\rho\le-2k$ in $D_{\rm ext}$.
Moreover, we get
\[
\begin{array}{l}
\d\left|\frac{\p}{\p\rho}\left(2\sqrt{2}\pi\sqrt{r}R(f|x,t)\right)\right|\\
\d=\left|\frac{\p}{\p\rho}\left(\int_{-5k/4}^{5k/4}
\frac{ds}{\sqrt{s-\rho+\rho^2z/2}}\int_{s=<\omega,y>-|y|^2z/2}f(y)dS_y\right)\right|\\
\\
\d=\left|\int_{-5k/4}^{5k/4}\frac{1-\rho z}{2(s-\rho+|\rho|^2z/2)^{3/2}}ds
\int_{s=<\omega,y>-|y|^2z/2}f(y)dS_y\right|\\
\d\le\frac{1+|\rho|z}{(|\rho|-5k/4)^{3/2}}\int_{-5k/4}^{5k/4}ds
\int_{s=<\omega,y>-|y|^2z/2}|f(y)|dS_y.
\end{array}
\]
Hence it follows from $|\rho|z=(t-r)/r\le 2$ for $(x,t)\in D_{ext}$ that
\[
\left|\frac{\p}{\p\rho}\left(2\sqrt{2}\pi\sqrt{r}R(f|x,t)\right)\right|
\le\frac{3}{(|\rho|-5k/4)^{3/2}}\|f\|_{L^1(\R^2)}.
\]
Summing up all the estimates and noticing that $\p/\p t=-\p/\p\rho$, we obtain
\[
\begin{array}{l}
\d\left|r^{1/2}u_L(x,t)
-\frac{1}{2\sqrt{2}\pi\sqrt{|\rho|+\rho^2z/2}}\int_{\R^2}g(x)dx\right|\\
\d\le
\frac{5k/4\|g\|_{L^1(\R^2)}+3\|f\|_{L^1(\R^2)}}{2\sqrt{2}\pi(|\rho|-5k/4)^{3/2}}
\end{array}
\]
in $D_{\rm ext}$.
It is trivial that (\ref{decay_est}) in $D_{\rm ext}$ follows from this inequality.
The proof is now complete.
\hfill$\Box$
\vskip10pt
\par
In what follows, we consider the following integral equations:
\begin{equation}
\label{IE}
u(x,t)=u^0(x,t)+L(F)(x,t)\quad\mbox{for}\ (x,t)\in\R^2\times[0,\infty),
\end{equation}
where we set $u^0:=\e u_L$ and
\begin{equation}
\label{L}
L(F)(x,t):=\frac{1}{2\pi}\int_{0}^{t}(t-\tau)\int_{|\xi|\le 1}\frac{F(x+(t-\tau)\xi,\tau)}{\sqrt{1-|\xi|^2}}d\xi d\tau
\end{equation}
for $F\in C(\R^2\times[0,\infty))$.
We note that $u$ in (\ref{L}) solves
\[
\left\{
\begin{array}{ll}
u_{tt}-\Delta u=F & \mbox{in}\ \R^2\times[0,\infty),\\
u(x,0)=\e f(x),\ u_t(x,0)=\e g(x), & x\in\R^2
\end{array}
\right.
\]
when $F\in C^2(\R^2\times[0,\infty))$.
Next two lemmas are useful to handle the radially symmetric functions.
\begin{lem}[John \cite{J55}]
\label{lm:Planewave}
Let $b\in C([0,\infty))$.
Then, the identity    
\begin{equation}
\label{Planewave}
\begin{array}{ll}
\d \int_{|\omega|=1}b(|x+\rho \omega|)dS_\omega
\d = 4\int_{|\rho-r|}^{\rho+r}\lambda h(\lambda, \rho ,r)
b(\lambda)d\lambda
\end{array} 
\end{equation}
holds for $x\in\R^2,\ r=|x|$ and $\rho>0$, where $h$ is defined by
\begin{equation}
\label{h}
h(\lambda, \rho ,r)
:=\{(\lambda+r)^2-\rho^2\}^{-1/2}\{\rho^2-(\lambda-r)^2\}^{-1/2}.
\end{equation}
\end{lem}
\par
See \cite{J55} for the proof of this lemma.
\begin{lem}[Agemi and Takamura \cite{AT92}]
\label{lem:int_new_L}
Let $L$ be a linear integral operator defined by {\rm (\ref{L})} 
and $\Psi=\Psi(|x|,t)\in C([0,\infty)^2),\ x\in\R^2$.
Then we have that
\begin{equation}
\label{int_new_L}
L\left(\Psi\right)(x,t)=L_{1}\left(\Psi\right)(r,t)+L_{2}\left(\Psi\right)(r,t),\ r=|x|,
\end{equation}
where $L_{i}\left(\Psi\right)\ (i=1,2)$ are defined by 
\begin{equation}
\label{int_domain1}
L_{1}\left(\Psi\right)(r,t):=
\d \frac{2}{\pi}\int_{0}^{t}d\tau\int_{|t-r-\tau|}^{t+r-\tau}
\hspace{-10pt}\lambda \Psi(\lambda,\tau)d\lambda
\int_{|\lambda-r|}^{t-\tau}
\frac{\rho h(\lambda,\rho,r)}{\sqrt{(t-\tau)^2-\rho^2}}d\rho,
\end{equation}
\begin{equation}
\label{int_domain2}
L_{2}\left(\Psi\right)(r,t):=
\frac{2}{\pi}\int_{0}^{(t-r)_{+}}\hspace{-20pt}d\tau
\int_{0}^{t-r-\tau}\hspace{-10pt}\lambda \Psi(\lambda,\tau)d\lambda
\int_{|\lambda-r|}^{\lambda+r}
\frac{\rho h(\lambda,\rho,r)}{\sqrt{(t-\tau)^2-\rho^2}}d\rho,
\end{equation}
where $a_+=\max\{a,0\}$.
Moreover, the following estimates hold in $[0,\infty)^2$:
\begin{equation}
\label{int_domain1_est}
|L_{1}\left(\Psi\right)(r,t)|\le 
\frac{1}{\sqrt{2r}}\int_{0}^{t}d\tau\int_{|r-t+\tau|}^{r+t-\tau}
\frac{\lambda |\Psi(\lambda,\tau)|d\lambda}{\sqrt{\tau+\lambda-t+r}},
\end{equation}
\begin{equation}
\label{int_domain2_est}
|L_{2}\left(\Psi\right)(r,t)|
 \le\int_{0}^{(t-r)_+}\hspace{-20pt}d\tau\int_{0}^{t-r-\tau}
\frac{\lambda |\Psi(\lambda,\tau)|d\lambda}
{\sqrt{t-r+\lambda-\tau}\sqrt{t-r-\tau-\lambda}}.
\end{equation}
\end{lem}
{\bf Proof.}
For the sake of completeness, we prove this lemma.
Changing variables by $y-x=(t-\tau)\xi$ in (\ref{L}),
we obtain that
\[
L\left(\Psi\right)(x,t)=\frac{1}{2\pi}\int_{0}^{t}d\tau\int_{|y-x|\le t-\tau}
\frac{\Psi(|y|,\tau)}{\sqrt{(t-\tau)^2-|y-x|^2}}dy.
\]
Introducing polar coordinates, we have that
\[
L\left(\Psi\right)(x,t)=\d \frac{1}{2\pi}\int_{0}^{t}d\tau\int_{0}^{t-\tau}
\frac{\rho d\rho}{\sqrt{(t-\tau)^2-\rho^2}}
\int_{|\omega|=1}
\Psi(|x+\rho\omega|,\tau)dS_{\omega}.
\]
Thus Lemma \ref{lm:Planewave} yields that
\begin{equation}
\label{int_new_L_last_even}
\begin{array}{llll}
L\left(\Psi\right)(x,t)& \d =\frac{2}{\pi}\int_{0}^{t}d\tau\int_{0}^{t-\tau}
\frac{\rho d\rho}{\sqrt{(t-\tau)^2-\rho^2}}\times&\\
&\quad \times \d \int_{|\rho-r|}^{\rho+r}\lambda 
\Psi(\lambda,\tau)h(\lambda,\rho,r)d\lambda.
\end{array}
\end{equation}
Therefore, (\ref{int_new_L}) follows from inverting the order of $(\rho,\lambda)$-integral
in (\ref{int_new_L_last_even}).
\par
The estimates (\ref{int_domain1_est}) and (\ref{int_domain2_est})
are established in the following way.
Note that
\begin{equation}
\label{L_1rho}
\lambda+r-\rho\ge \lambda+r-t+\tau,
\quad \lambda+r+\rho\ge \lambda+r+|\lambda-r|\ge 2r
\end{equation}
for $\rho\le t-\tau$ and $\rho\ge |\lambda-r|$.
It is easy to see that
\begin{equation}
\label{beta_func}
\int_{a}^{b}\frac{\rho d\rho}{\sqrt{\rho^2-a^2}\sqrt{b^2-\rho^2}}=\frac{\pi}{2}\ \mbox{for}\ 0\le a <b. 
\end{equation}
Hence (\ref{int_domain1_est}) follows from
(\ref{L_1rho}) and (\ref{beta_func}) with $a=|\lambda-r|$ and $b=t-\tau$.
Next let $t>r$. Since we have that
\begin{equation}
\label{L_2rho1}
t-\tau-\rho\ge t-\tau-\lambda-r
\end{equation}
for $\rho\le \lambda+r$ and that
\begin{equation}
\label{L_2rho2}
t-\tau+\rho\ge t-\tau+|\lambda-r|\ge t-r-\tau+\lambda 
\end{equation}
for $\rho\ge |\lambda-r|$, $\lambda\le t-r-\tau$,
we obtain (\ref{int_domain2_est}) by
(\ref{L_2rho1}), (\ref{L_2rho2}) and (\ref{beta_func}) with $a=|\lambda-r|$ and $b=\lambda+r$.
\hfill$\Box$

\section{A priori estimate}
In this section, we show a priori estimates
which play key roles in the classical iteration method as in John \cite{J79}. 
First of all, we define some weighted $L^{\infty}$ norms.

For $r,t\ge0$, we define the following weighted functions:
\begin{eqnarray}
\label{weight1}
w_1(r,t)&:=&\tau_+(r,t)^{1/2}\tau_-(r,t)^{1/2},\\
\label{weight2}
w_2(r,t)&:=&\tau_+(r,t)^{1/2}\tau_-(r,t)^{3/2},\\
\label{weight3}
w_3(r,t)&:=&\tau_+(r,t)^{p/2-1},
\end{eqnarray}
where we set
\[
\tau_+(r,t):=\frac{t+r+2k}{k}, \quad\tau_-(r,t):=\frac{t-r+2k}{k}.
\]
For these weighted functions, we denote weighted $L^\infty$ norms of $V$ by 
\begin{equation}
\label{norm}
\|V\|_i:=\sup_{(x,t)\in\R^2\times[0,T]}\{w_i(|x|,t)|V(x,t)|\},
\end{equation}
where $i=1,2,3$.

The following lemma is one of the most essential estimates.
\begin{lem}
\label{lm:apriori}
Let $L$ be the linear integral operator defined by {\rm (\ref{L})}.
Assume that $V\in C(\R^2\times[0,T])$ with
{\rm  supp} $V\subset\{(x,t)\in \R^2\times[0,T] : |x|\le t+k\}$ and $\|V\|_i<\infty$ $(i=1,3)$.
Then, there exists a positive constant $C_1$ independent of $k$ and $T$ such that 
\begin{equation}
\label{apriori1}
\|L(|V|^p)\|_1\le C_1k^2\|V\|_1^pD_1(T),
\end{equation}
\begin{equation}
\label{apriori2}
\|L(|V|^p)\|_3\le C_1k^2\|V\|_3^pD_2(T),
\end{equation}
where $D_i(T)$ $(i=1,2)$ are defined by 
\begin{equation}
\label{D_1}
D_1(T):=
\left(\frac{2T+3k}{k}\right)^{3-p},\\
\end{equation}
\begin{equation}
\label{D_2}
D_2(T):=\left(\frac{2T+3k}{k}\right)^{\gamma(p,2)/2}.
\end{equation}
\end{lem}
{\bf Proof of  Lemma \ref{lm:apriori}.}
The proof is divided into two pieces according  to (\ref{apriori1}) and (\ref{apriori2}). 
From now on, a positive constant $C$ independent of $\e$ and $k$
may change from line to line.
\par\noindent
{\bf Estimate in (\ref{apriori1}).}
It is clear that (\ref{apriori1}) follows from the basic estimate:
\begin{equation}
\label{basic1}
L(w_1^{-p})\le C k^2 w_1(r,t)^{-1} D_1(T).
\end{equation}
First we shall show a part of (\ref{basic1}),
\[
L_1(w_1^{-p})\le C k^2 w_1(r,t)^{-1} D_1(T).
\]
Introducing the characteristic variables
\begin{equation}
\label{alpha-beta}
\alpha=\tau+\lambda,\quad \beta=\tau-\lambda
\end{equation}
in the integral of (\ref{int_domain1_est}), we get
\[
\begin{array}{lll}
|L_1(\Psi)|
&\d \le \frac{C}{\sqrt{r}}
\int_{|t-r|}^{t+r}d\alpha 
\int_{-k}^{t-r}\left|\Psi\left(\frac{\alpha-\beta}{2}, \frac{\alpha+\beta}{2}\right)\right|
\frac{(\alpha-\beta)}{\sqrt{\alpha-t+r}}d\beta.&
\end{array}
\]
Setting $\Psi(\lambda,\tau)=\{w_1(\lambda,\tau)\}^{-p}$, we have
\[
\begin{array}{lll}
L_1(w_1^{-p})
&\d \le \frac{Ck}{\sqrt{r}}\int_{t-r}^{t+r}\left(\frac{\alpha+2k}{k}\right)^{1-p/2}
\hspace{-20pt}\frac{d\alpha}{\sqrt{\alpha-t+r}}
\int_{-k}^{t-r}\left(\frac{\beta+2k}{k}\right)^{-p/2}d\beta&\\
&\d \le Ck^2\tau_+(r,t)^{1-p/2}\tau_+(r,t)^{1-p/2}.&
\end{array}
\]
Thus, recalling (\ref{weight1}) and (\ref{D_1}), we get
\[
L_1(w_1^{-p})\le Ck^2w_1(r,t)^{-1}\tau_+(r,t)^{3-p}\le Ck^2w_1(r,t)^{-1}D_1(T).
\]
\par
Next we shall show the remaining part of (\ref{basic1}),
\[
L_2(w_1^{-p})\le C k^2w_1(r,t)^{-1}D_1(T)
\quad \mbox{for}\ t-r\ge0.
\]
Introducing the characteristic variables (\ref{alpha-beta}) in the integral of (\ref{int_domain2_est}), we get
\[
|L_2(\Psi)|
\d \le C
\int_{0}^{t-r}d\alpha\int_{-k}^{t-r}
\left|\Psi\left(\frac{\alpha-\beta}{2}, \frac{\alpha+\beta}{2}\right)\right|
\frac{(\alpha-\beta)}{\sqrt{t-r-\alpha}}\frac{d\beta}{\sqrt{t-r-\beta}}.
\]
Setting $\Psi(\lambda,\tau)=\{w_1(\lambda,\tau)\}^{-p}$, we have that
\begin{equation}
\label{L_2est1}
\begin{array}{lll}
L_2(w_1^{-p})
&\d \le Ck\int_{0}^{t-r}\left(\frac{\alpha+2k}{k}\right)^{1-p/2}\frac{d\alpha}{\sqrt{t-r-\alpha}}\times\\
&\quad \d \times\int_{-k}^{t-r}\left(\frac{\beta+2k}{k}\right)^{-p/2}\frac{d\beta}{\sqrt{t-r-\beta}}.
\end{array}
\end{equation}
First we consider the case of $t-r\ge k$. Then, we get
\begin{equation}
\label{J_1-J_2}
L_2(w_1^{-p})\le Ck\tau_{+}(r,t)^{1-p/2}\{J_1(r,t)+J_2(r,t)\},
\end{equation}
where we set
\[
J_1(r,t):=\int_{0}^{t-r}\frac{d\alpha}{\sqrt{t-r-\alpha}}
\int_{-k}^{(t-r-k)/2}\left(\frac{\beta+2k}{k}\right)^{-p/2}\frac{d\beta}{\sqrt{t-r-\beta}},
\]
\[
J_2(r,t):=\int_{0}^{t-r}\frac{d\alpha}{\sqrt{t-r-\alpha}}
\int_{(t-r-k)/2}^{t-r}\left(\frac{\beta+2k}{k}\right)^{-p/2}\frac{d\beta}{\sqrt{t-r-\beta}}.
\]
It is easy to see that
\begin{equation}
\label{J_1est}
\begin{array}{ll}
\d 
J_1(r,t)
&\d \le \frac{C}{\sqrt{t-r+k}}\int_{0}^{t-r}\hspace{-10pt}\frac{d\alpha}{\sqrt{t-r-\alpha}}
\int_{-k}^{t-r}\left(\frac{\beta+2k}{k}\right)^{-p/2}\hspace{-10pt}d\beta\\
&\\
&\d \le Ck\tau_{+}(r,t)^{1-p/2}.
\end{array}
\end{equation}
On the other hand, it follows from $1-p/2>0$ that
\begin{equation}
\label{J_2est}
\begin{array}{lll}
J_2(r,t)
&\d \le C\tau_{-}(r,t)^{-p/2}
\int_{0}^{t-r}\frac{d\alpha}{\sqrt{t-r-\alpha}}
\int_{(t-r-k)/2}^{t-r}\frac{d\beta}{\sqrt{t-r-\beta}}\\
&\d \le Ck\tau_{-}(r,t)^{1-p/2}\le Ck \tau_{+}(r,t)^{1-p/2}.
\end{array}
\end{equation}
Hence (\ref{J_1-J_2}), (\ref{J_1est}) and (\ref{J_2est}) yield that
\[
L_2(w_1^{-p})\le Ck^2\tau_{+}(r,t)^{2-p}\le Ck^2 w_1(r,t)^{-1}D_1(T)
\]
for $t-r\ge k$.
In the case of $0<t-r\le k$, (\ref{L_2est1}) yields that
\[
\begin{array}{lll}
L_2(w_1^{-p})
&\d \le Ck\int_{0}^{t-r}\left(\frac{\alpha+2k}{k}\right)^{-p/2}\frac{d\alpha}{\sqrt{t-r-\alpha}}\times\\
&\quad \d \times\int_{-k}^{t-r}\left(\frac{\beta+2k}{k}\right)^{-p/2}\frac{d\beta}{\sqrt{t-r-\beta}}\\
&\d \le Ck^2=Ck^2 w_1(r,t)^{-1}w_1(r,t)\le Ck^2 w_1(r,t)^{-1} D_1(T).
\end{array}
\]
Therefore (\ref{basic1}) is now established.
\par\noindent
{\bf Estimate in  (\ref{apriori2}).}
Similarly to the above, we note that
(\ref{apriori2}) follows from the basic estimate:
\begin{equation}
\label{basic2}
L(w_3^{-p})\le C k^2 w_3(r,t)^{-1} D_2(T).
\end{equation}
Setting $\Psi(\lambda,\tau)=\{w_3(\lambda,\tau)\}^{-p}$ in (\ref{int_domain1_est}) and 
introducing (\ref{alpha-beta}) in the integral of (\ref{int_domain1_est}), we get
\[
\begin{array}{ll}
L_1(w_3^{-p})
&\d\le \frac{Ck}{\sqrt{r}}
\int_{|t-r|}^{t+r}\left(\frac{\alpha+2k}{k}\right)^{1+(1-p/2)p}\frac{d\alpha}{\sqrt{\alpha-t+r}}
\int_{-k}^{t-r}d\beta\\
&\d \le Ck^2 \tau_+(r,t)^{(1-p/2)p+1}\tau_{-}(r,t)
\le Ck^2\tau_+(r,t)^{1-p/2+\gamma(2,p)/2}.
\end{array}
\]
Hence we have a part of (\ref{basic2}),
\[
L_1(w_3^{-p})\le C k^2 w_3(r,t)^{-1}D_2(T).
\]
Now, let $t-r>0$.
Setting $\Psi(\lambda,\tau)=\{w_3(\lambda,\tau)\}^{-p}$ in (\ref{int_domain2_est}) and 
introducing (\ref{alpha-beta}) in the integral of (\ref{int_domain2_est}), we get
\[
\begin{array}{lll}
L_2(w_3^{-p})
&\d \le Ck\int_{0}^{t-r}\left(\frac{\alpha+2k}{k}\right)^{1+(1-p/2)p}
\frac{d\alpha}{\sqrt{t-r-\alpha}}\int_{-k}^{t-r}\frac{d\beta}{\sqrt{t-r-\beta}}\\
&\d\le Ck\tau_{+}(r,t)^{(1-p/2)p+1}
\int_{0}^{t-r}\frac{d\alpha}{\sqrt{t-r-\alpha}}
\int_{-k}^{t-r}\frac{d\beta}{\sqrt{t-r-\beta}}\\
&\d \le Ck^2\tau_{+}(r,t)^{(1-p/2)p+2}.
\end{array}
\]
Therefore we have that
\[
L_2(w_3^{-p})\le C k^2w_3(r,t)^{-1}D_2(T)
\]
for $t-r\ge0$ which yields (\ref{basic2}).
The proof of Lemma \ref{lm:apriori} is complete.
\hfill$\Box$
\vskip10pt
\par
In order to construct a solution in our weighted $L^\infty$ space,
the following variant to the a priori estimate is required.

\begin{lem}
\label{lm:apriori_u0}
Let $L$ be the linear integral operator defined by {\rm (\ref{L})} and $0\le \nu<p$.
Assume that $V,V_0 \in C(\R^2\times[0,T])$
with {\rm  supp} $(V,V_0)\subset\{(x,t)\in \R^2\times[0,T] : |x|\le t+k\}$,
and $\|V_{0}\|_2, \|V\|_3<\infty$. 
Then, there exists a positive constant $C_2$ independent of $k$ such that 
\begin{equation}
\label{apriori_u0}
\|L(|V_0|^{p-\nu}|V|^{\nu})\|_3
\le C_2k^2\|V_0\|_2^{p-\nu}\|V\|_3^{\nu}D_{3,\nu}(T),
\end{equation}
where $D_{3,\nu}(T)$ are defined by 
\begin{equation}
\label{D_3}
D_{3,\nu}(T):=
\left\{
\begin{array}{lll}
\d \left(\frac{2T+3k}{k}\right)^{\nu(3-p)/2} & \mbox{if} &\d p>\nu+\frac{2}{3},\\
\d \log\frac{2T+3k}{k}\left(\frac{2T+3k}{k}\right)^{(7/3-\nu)\nu/2} & \mbox{if} &\d p=\nu+\frac{2}{3},\\
\d \left(\frac{2T+3k}{k}\right)^{1-3p/2+(3-p/2)\nu} & \mbox{if} &\d p<\frac{2}{3}+\nu.
\end{array}
\right.
\end{equation}
\end{lem}
{\bf Proof.}
Similarly to the proof of Lemma \ref{lm:apriori}, we note that
(\ref{apriori_u0}) follows from
\begin{equation}
\label{basic3_1}
L(w_2^{-(p-\nu)}w_3^{-\nu})\le C k^2 w_3(r,t)^{-1}D_{3,\nu}(T).
\end{equation}
A part of this estimate is
\begin{equation}
\label{L_1basic1}
L_1(\tau_{+}^{-(p-\nu)/2+(1-p/2)\nu}\tau_{-}^{-3(p-\nu)/2})\le C k^2 w_3(r,t)^{-1}D_{3,\nu}(T).
\end{equation}
Setting 
\begin{equation}
\label{weight-set1}
\Psi(\lambda,\tau)=\tau_{+}^{-(p-\nu)/2+(1-p/2)\nu}\tau_{-}^{-3(p-\nu)/2}
\end{equation}
in (\ref{int_domain1_est}) and introducing (\ref{alpha-beta}) in the integral in (\ref{int_domain1_est}), we get
\begin{equation}
\label{L_1basic1_est1}
\begin{array}{l}
L_1(\tau_{+}^{-(p-\nu)/2+(1-p/2)\nu}\tau_{-}^{-3(p-\nu)/2})\\
\\
\d\le \frac{Ck}{\sqrt{r}}
\int_{|t-r|}^{t+r}\left(\frac{\alpha+2k}{k}\right)^{1-(p-\nu)/2+(1-p/2)\nu}\frac{d\alpha}{\sqrt{\alpha-t+r}}\\
\d \quad\times\int_{-k}^{t-r}\left(\frac{\beta+2k}{k}\right)^{-3(p-\nu)/2}d\beta.
\end{array}
\end{equation}
Note that 
\begin{equation}
\label{L_1_beta-est1}
\d\int_{-k}^{t-r}\left(\frac{\beta+2k}{k}\right)^{-3(p-\nu)/2}\hspace{-20pt}d\beta\le
\left\{
\begin{array}{lll}
Ck & \mbox{if} &\d p>\nu+2/3,\\
\d Ck\log \tau_{-}(r,t)& \mbox{if} &\d p=\nu+2/3,\\
\d Ck\tau_{-}(r,t)^{1-3(p-\nu)/2} & \mbox{if} &\d p<\nu+2/3
\end{array}
\right.
\end{equation}
and that the $\alpha$-integral in (\ref{L_1basic1_est1})
 is estimated by
 \[
 Ck \sqrt{r}\tau_+(r,t)^{1-(p-\nu)/2+(1-p/2)\nu}.
 \]
Therefore   (\ref{L_1basic1}) is follows from $\tau_+^{1-(p-\nu)/2+(1-p/2)\nu}\le w_3^{-1}\tau_+^{(3-p)\nu/2}$.
\par
Next, we shall show the remaining part of (\ref{basic3_1}),
\begin{equation}
\label{L_2basic1}
L_2(\tau_{+}^{-(p-\nu)/2+(1-p/2)\nu}\tau_{-}^{-3(p-\nu)/2})
\le C k^2 w_3(r,t)^{-1}D_{3,\nu}(T).
\end{equation}
Setting $\Psi(\lambda,\tau)$ such as (\ref{weight-set1}) in (\ref{int_domain2_est}) and 
introducing (\ref{alpha-beta}) in the integral of (\ref{int_domain2_est}), we get
\begin{equation}
\label{L_2basic1_est1}
\begin{array}{l}
L_2(\tau_{+}^{-(p-\nu)/2+(1-p/2)\nu}\tau_{-}^{-3(p-\nu)/2})\\
\\
\d \le Ck\int_{0}^{t-r}\left(\frac{\alpha+2k}{k}\right)^{1-(p-\nu)/2+(1-p/2)\nu}
\frac{d\alpha}{\sqrt{t-r-\alpha}}\\
\quad \d \times\int_{-k}^{t-r}\left(\frac{\beta+2k}{k}\right)^{-3(p-\nu)/2}\frac{d\beta}{\sqrt{t-r-\beta}}.
\end{array}
\end{equation}
First we consider the case of $t-r\ge k$. Then we have
\[
\begin{array}{lll}
L_2(\tau_{+}^{-(p-\nu)/2+(1-p/2)\nu}\tau_{-}^{-3(p-\nu)/2})\\
\le Ckw_3(r,t)^{-1}\tau_{+}(r,t)^{(3-p)\nu/2}\{K_1(r,t)+K_2(r,t)\},
\end{array}
\]
where we set
\[
K_1(r,t):=\int_{0}^{t-r}\frac{d\alpha}{\sqrt{t-r-\alpha}}
\int_{-k}^{(t-r-k)/2}\left(\frac{\beta+2k}{k}\right)^{-3(p-\nu)/2}\frac{d\beta}{\sqrt{t-r-\beta}},
\]
\[
K_2(r,t):=\int_{0}^{t-r}\frac{d\alpha}{\sqrt{t-r-\alpha}}
\int_{(t-r-k)/2}^{t-r}\left(\frac{\beta+2k}{k}\right)^{-3(p-\nu)/2}\frac{d\beta}{\sqrt{t-r-\beta}}.
\]
It follows from (\ref{L_1_beta-est1}) that
\[
K_1(r,t)
\le C(t-r+k)^{-1/2}\int_{0}^{t-r}\hspace{-15pt}\frac{d\alpha}{\sqrt{t-r-\alpha}}
\int_{-k}^{t-r}\left(\frac{\beta+2k}{k}\right)^{-3(p-\nu)/2}d\beta.
\]
Hence (\ref{L_2basic1}) for $t-r\ge k$ is established by (\ref{L_1_beta-est1}).
Moreover it is easy to see that
\[
\begin{array}{ll}
K_2(r,t)
&\d \le C\tau_{-}(r,t)^{-3(p-\nu)/2}
\int_{0}^{t-r}\frac{d\alpha}{\sqrt{t-r-\alpha}}
\int_{(t-r-k)/2}^{t-r}\frac{d\beta}{\sqrt{t-r-\beta}}\\
&\d \le Ck\tau_{-}(r,t)^{1-3(p-\nu)/2}.
\end{array}
\]
Hence (\ref{L_2basic1}) for $t-r\ge k$ is established.
In the case of $0\le t-r\le k$, (\ref{L_2basic1_est1}) yields that
\[
\begin{array}{l}
L_2(\tau_{+}^{-(p-\nu)/2+(1-p/2)\nu}\tau_{-}^{-3(p-\nu)/2})\\
\\
\d \le Ck\int_{0}^{t-r}\left(\frac{\alpha+2k}{k}\right)^{1-(p-\nu)/2+(1-p/2)\nu}
\frac{d\alpha}{\sqrt{t-r-\alpha}}\\
\quad \d \times\int_{-k}^{t-r}\left(\frac{\beta+2k}{k}\right)^{-3(p-\nu)/2}\frac{d\beta}{\sqrt{t-r-\beta}}\\
\\
\d \le Ck^2=Ck^2 w_3(r,t)^{-1}w_3(r,t)\le Ck^2 w_3(r,t)^{-1}D_{3,\nu}(T).
\end{array}
\]
This gives us (\ref{L_2basic1}) for $0\le t-r\le k$
which leads to (\ref{basic3_1}).
The proof of Lemma \ref{lm:apriori_u0} is complete.
\hfill$\Box$
\vskip10pt
\par
In order to pick up the sufficient condition to construct a solution
in our weighted $L^\infty$ space, we need the following  lemma
on comparison among the quantities depending on $T$.

\begin{lem}
\label{lm:D_2-D_3}
Let $1<p<2$ and let $D_2(T)$ and $D_{3,\nu}(T)$ are the one in (\ref{D_2}) and (\ref{D_3}). Then we have
\begin{equation}
\label{D_2-D_3_1}
D_{3,1}(T)\le D_2(T)^{1/p},
\end{equation}
\begin{equation}
\label{D_2-D_3_2}
D_{3,p-1}(T)\le D_2(T)^{(p-1)/(p+1)},
\end{equation}
\begin{equation}
\label{D_2-D_3_3}
D_{3,0}(T)=1.
\end{equation}
\end{lem}
\par\noindent
{\bf Proof.}
All the estimates follow from direct computations.
When $p>5/3$, we have that $(3-p)p/2<\gamma(p,2)/2$
which yields (\ref{D_2-D_3_1}). 
When $p=5/3$, we have that $(7/3-\nu)\nu/2=(3-p)\nu/2$ for $p=\nu+2/3$, 
so that it is sufficient to show that $p\{\delta+(3-p)/2\}<\gamma(p,2)/2$ holds for suitable $\delta>0$ to get  (\ref{D_2-D_3_1}). 
This can be guaranteed by taking $0<\delta<1/p$.
When $1<p<5/3$,  we have that $p(4-2p)<\gamma(p,2)/2$ holds.
Therefore (\ref{D_2-D_3_1}) is established in all the cases.
\par
The estimate (\ref{D_2-D_3_2}) follows from $(p+1)(3-p)/2<\gamma(p,2)/2$
which is equivalent to $p>1$.
The estimate (\ref{D_2-D_3_3}) is trivial by definition of $D_{3,\nu}(T)$.
The proof is now complete. 
\hfill $\Box$

\section{Lower bound of the lifespan}
First, we shall show the estimate for (\ref{lower_lifespan}) in the case of
\newline
$\d\int_{\R^2} g(x) dx=0$. We consider
the following integral equation:
\begin{align}
U=L(|u^{0}+U|^{p}) \quad \mbox{in} \quad \R^{2} \times [0, T]
\label{dai1}
\end{align}
where $L$ and $u^{0}$ are defined in (\ref{IE}). 
Suppose we obtain the solution $U(x,t)$ of (\ref{dai1}).
Then, putting $u=U+u^0$, we get the solution of (\ref{IE})
and its life span is the same as that of $U$. Thus  we have
reduced the problem to the analysis of (\ref{dai1}).
In view of (\ref{IE}) and (\ref{L}), we note that $\partial U/\partial t$ can be expressed in
$\nabla_{x} U$. Hence we consider spatial derivatives of $U$ only.

We define $U_{l}$ by
\begin{align}
	U_{1}=0,\quad U_{l}=L(|u^0+U_{l-1}|^p)\quad \mbox{for} \quad 
	l \geq 2.\label{dai2}
\end{align}
We take $\e$ and $T$ such that
\begin{align}
	C\e^{p(p-1)}D_2(T) \leq 1\label{dai3}
\end{align}
where
\begin{align}
	&C:=(2^{2p}p)^{\frac{p}{p-1}}
	\max \{C_{1}k^2M_{0}^{p-1}, (C_{2}k^2C_{0}^{p-1})^{p},
	(C_{2}k^2M_{0}^{p-2}C_{0})^{\frac{p}{p-1}} \},\label{dai4}\\
	&M_{0}:=2^{p} p C_{3} k^2 C_{0}^p\label{dai6}
\end{align}
and $C_{3}:=\max\{C_{1},C_{2}\}$.

In order to get a $C^{1}$ solution of (\ref{dai1}), we shall show the
convergence of $ \{ U_{l} \}_{ l \in \N }$ in a function space
$X$ defined by
\begin{align}
	X:=\{& U(x,t) \ : \ \nabla_x^{\alpha}U \in C(\R^{2} \times [0,T]) \ 
	\mbox{for} \ |\alpha| \leq 1,\
	 \|  U \|_{X} < \infty,\nonumber\\ 
	 &\mbox{supp} \ U
	\subset \{(x,t):|x| \leq t+k \} \}\nonumber
\end{align}
which is equipped with the norm
\begin{align}
	\| U \|_{X} =
	 \sum_{|\alpha| \leq 1} \| \nabla_{x}^{\alpha} U \|_3.\nonumber
\end{align}
We see that $X$ is a Banach space for any fixed $T > 0$.
It follows from the definition of the norm (\ref{norm})
that there exists a positive constant
$C_{T}$ depending on $T$ such that
\begin{align}
	\| U \|_{3} \geq C_{T} | U(x,t)|,\quad t \in [0,T].\nonumber
\end{align}

By induction, we shall obtain
\begin{align}
	\| U_{l} \|_{3} \leq 2M_{0}  \e^p.\label{dai5}
\end{align}
For $l=1$, (\ref{dai5}) holds. Assume that $\| U_{l-1}\|_3 \leq 2M_{0} \e^{p}$ $(l \geq 2)$.
Since 
\[
|a+b|^p\le 2^{p-1}(|a|^p+|b|^p)\quad\mbox{for}\ p>1\ \mbox{and }\ a,b\in \R,
\]
we get from (\ref{dai2}) that
\begin{align}
	\| U_{l} \|_3\le 2^{p-1}\{\|L(|u^0|^p)\|_3
	+\|L(|U_{l-1}|^p)\|_3\}\label{dai7}.
\end{align}
Making use of (\ref{apriori_u0}) with $\nu=0$, (\ref{D_2-D_3_3}) and (\ref{decay_est}),
we have
\begin{align}
	\|L(|u^0|^p)\|_3 &\le C_2k^2\|u^0\|_2^p\nonumber\\
	&\le C_2 k^2C_0^p\e^p,\label{dai8}
\end{align}
where we used (\ref{decay_est}) with $\d\int_{\R^2}g(x)dx=0$,
(\ref{weight2}) and (\ref{norm}).
We see from (\ref{apriori2}) that
\begin{align}
	\|L(|U_{l-1}|^p)\|_3 &\le 
	C_1k^2\|U_{l-1}\|_3^pD_2(T)\nonumber\\
	&\le C_1k^2(2M_{0} \e^p)^pD_2(T).\label{dai9}
\end{align}
Summarizing (\ref{dai7}), (\ref{dai8}), (\ref{dai9})
and (\ref{dai6}), we get
\begin{align}
	\| U_{l} \|_{3} \le M_{0} \e^{p}
	+2^{p-1}C_1k^2(2M_{0} \e^p)^pD_2(T).\label{dai10}
\end{align}
This inequality shows (\ref{dai5}) provided
(\ref{dai3}) and (\ref{dai4}) hold.

We shall estimate the differences of $\{ U_{l} \}_{l \in \N}$.
Since 
\[
||a|^p-|b|^p| \leq p(|a|^{p-1}+|b|^{p-1})|a-b|\quad\mbox{for}\ p>1,
\]
we obtain from (\ref{dai2}) that
\begin{align}
	\| U_{l+1} -U_{l} \|_{3} &\leq 2^{p-1} p
	\{2\| L(|u^{0}|^{p-1} | U_{l}-U_{l-1}| )\|_{3}\nonumber\\
	&+\| L((|U_{l}|^{p-1}+|U_{l-1}|^{p-1}) 
	|U_{l}-U_{l-1}|)\|_{3}\}.\label{dai12}
\end{align}
From (\ref{apriori_u0}) with $\nu=1$, (\ref{D_2-D_3_1})
and (\ref{decay_est}), we obtain
\begin{align}
	\| L(|u^0|^{p-1} |U_{l}-U_{l-1}|)\|_{3}
	&\leq C_{2}k^2\| u^{0} \|_{2}^{p-1} \| U_{l} - U_{l-1} \|_{3}
	D_{2}(T)^{1/p}\notag\\
	&\leq C_{2}k^2 (C_0 \e)^{p-1} D_{2}(T)^{1/p}
	\| U_{l} -U_{l-1} \|_{3}.\label{dai14}
\end{align}
We get from (\ref{apriori2}) and (\ref{dai5}) that
\begin{align}
	&\| L((|U_{l}|^{p-1} +|U_{l-1}|^{p-1})
	|U_{l}-U_{l-1}|)\|_{3}\nonumber\\
	\leq &C_{1}k^2D_{2}(T) 
	(\|  U_{l} \|_{3}^{p-1}+\|  U_{l} \|_{3}^{p-1})
	\| U_{l}-U_{l-1} \|_{3}\notag\\
	\leq &2C_{1}k^2D_{2}(T)(2M_{0} \e^p)^{p-1}
	\| U_{l} - U_{l-1} \|_{3}.\label{dai15}
\end{align}
Hence, we obtain from (\ref{dai12}), (\ref{dai14}) and
(\ref{dai15}) that
\[
	\| U_{l+1} - U_{l} \|_{3} 
	\leq \frac{1}{2}\| U_{l}-U_{l-1}\|_{3}
\]
provided (\ref{dai3}) and (\ref{dai4}) hold.
Therefore we have
\begin{align}
	\| U_{l+1}-U_{l} \|_{3} \leq 2^{-l}C_{4}\quad 
	\mbox{for} \quad l\geq1.\label{dai18}
\end{align}
\par
Next, by induction, we shall show the following boundedness
of $\{ \partial_i U_l\}$.
\begin{align}
	\| \partial_{i} U_{l}\|_{3} \leq 2M_{0} \e^p.\label{dai19}
\end{align}
Assume that $\| \partial_i U_{l-1} \|_{3} \leq 2M_{0} \e^p$ \ 
$(l \geq 2)$. From (\ref{dai2}), we have
\begin{align}
	\| \partial_{i} U_{l} \| _3\leq 2^{p-1}p
	\| L\{ (|u^{0}|^{p-1}+|U_{l-1}|^{p-1})
	(|\partial_{i} u^{0}|+|\partial_{i} U_{l-1}|)\} \|_{3}.\label{dai20}
\end{align}
By using Lemma \ref{lm:apriori_u0} and 
Lemma \ref{lm:D_2-D_3}, we shall show
\begin{align}
	&\| L(|u^0|^{p-1}|\partial_{i}u^0|)\|_{3}
	\leq C_{2}k^2(C_{0} \e)^p,\label{dai21}\\
	&\| L(|U_{l-1}|^{p-1}|\partial_{i}u^0|)\|_{3}
	\leq C_{2}k^2(2M_{0} \e^p)^{p-1}
	C_{0} \e D_2(T)^{\frac{p-1}{p+1}},\label{dai22}\\
	&\| L(|u^0|^{p-1}|\partial_{i}U_{l-1}|)\|_{3}
	\leq C_{2}k^2(C_{0} \e)^{p-1}D_{2}(T)^{1/p}
	(2M_{0} \e^p).\label{dai23}
\end{align}
We shall prove only (\ref{dai21}), since we can
prove (\ref{dai22}) and (\ref{dai23}) in a similar way.
It follows from  (\ref{apriori_u0}) with $\nu=0$, 
(\ref{D_2-D_3_3}) and (\ref{decay_est}) that
\begin{align}
	\| L(|u^0|^{p-1}|\partial_{i}u^0|)\|_{3}
	&\leq C_{2}k^2
	\| u^0\|_{2}^{p-1} \| \partial_{i} u^0\|_{2}\nonumber\\
	&\leq C_{2}k^2(C_{0} \e)^p.\nonumber
\end{align}
We have from (\ref{apriori2}) and (\ref{dai5}) that
\begin{align}
	\| L(|U_{l-1}|^{p-1}|\partial_{i}U_{l-1}|)\|_{3}
	&\leq C_{1}k^2\| U_{l-1}\|_{3}^{p-1} D_{2}(T)
	\| \partial_{i} U_{l-1}\|_{3}\nonumber\\
	&\leq C_{1}k^2(2 M_{0} \e^p)^{p}D_{2}(T).\label{dai24}
\end{align}
Summarizing (\ref{dai20}), (\ref{dai21}), (\ref{dai22}),
(\ref{dai23}) and (\ref{dai24}), we get from
\newline $D_{2}(T) \geq 1$ that
\begin{align}
	\| \partial_{i} U_l \|_{3}
	&\leq M_{0}\e^p
	+2^{p-1}pC_{3}k^2
	\{(2M_{0})^{p-1}C_{0}\epsilon^{p^2-p+1}D_2(T)^{\frac{p-1}{p}}\nonumber\\
	&+2M_{0}C_{0}^{p-1}\e^{2p-1}D_{2}(T)^{1/p}
	+(2M_{0})^p\e^{p^2}D_{2}(T)
	\}.\label{dai25}
\end{align}
Therefore we obtain (\ref{dai19}) provided (\ref{dai3}) and (\ref{dai4}) hold.
\par
Finally, we shall estimate the difference of ${\partial_{i}U_l}$.
We obtain from (\ref{dai2}) that
\begin{align}
	\| \partial_{i}(U_{l+1}-U_{l})\|_{3}
	&\leq  p\{
	2^{p-1}\| L((|u^0|^{p-1}+|U_{l}|^{p-1})
	|\partial_{i}(U_{l}-U_{l-1})|) \|_{3}\nonumber\\
	&+  \| L(|U_{l}-U_{l-1}|^{p-1}|\partial_{i}(u_{0}+U_{l-1})|) \|_{3}
	\}.\label{dai28}
\end{align}
We get from (\ref{apriori2}), (\ref{apriori_u0}) with $\nu=1$
and (\ref{D_2-D_3_1}) that
\begin{align}
	& 
	\| L((|u^0|^{p-1}+|U_{l}|^{p-1})
	|\partial_{i}(U_{l}-U_{l-1})|) \|_{3}\nonumber\\
	\leq & C_{3}k^2(\| u^{0} \|_{2}^{p-1}D_{2}(T)^{1/p}+
	\| U_{l} \|_3^{p-1}D_{2}(T))
	\| \partial_{i} (U_{l}-U_{l-1})\|_{3}\nonumber\\
	\leq & 
	C_{3}k^2\{ (C_{0} \e)^{(p-1)}D_{2}(T)^{1/p}
	+(2M_{0}\e^p)^{p-1}D_{2}(T)\}
	\| \partial_{i}(U_{l}-U_{l-1}) \|_3.\label{dai29}
\end{align}
It follows from (\ref{apriori2}), (\ref{apriori_u0}) with $\nu=p-1$
and (\ref{D_2-D_3_2}) that
\begin{align}
	& \| L(|U_{l}-U_{l-1}|^{p-1}
	|\partial_{i}(u^0+U_{l-1})|) \|_{3}\notag\\
	\leq & C_{3}k^2\| U_{l}-U_{l-1} \|_{3}^{p-1}
	\{
	\| \partial_{i} u^{0} \|_{2}D_{2}(T)^{\frac{p-1}{p+1}}
	+\| \partial_{i} U_{l-1} \|_{3} D_{2}(T)
	\}\notag\\
	\leq & C_{3}k^2\{
	C_{0} \e D_{2}(T)^{\frac{p-1}{p+1}}
	+2M_{0} \e^{p}D_{2}(T)	\}
	\| U_{l}-U_{l-1} \|_{3}^{p-1}\notag\\
	\leq & C_{5}\| U_{l}-U_{l-1} \|_{3}^{p-1}.\label{dai30}
\end{align}
Summarizing (\ref{dai28}), (\ref{dai29}), (\ref{dai30})
(\ref{dai4}) and (\ref{dai18}), we get
\begin{align}
	\| \partial_{i}(U_{l+1}-U_{l}) \|_{3}
	&\leq pC_{5}\| U_{l}-U_{l-1} \|_{3}^{p-1}
	+\frac{1}{2}\| \partial_{i}(U_{l}-U_{l-1}) \|_{3}\notag\\
	&\leq C_{6}2^{-l(p-1)}+\frac{1}{2}
	\| \partial_{i}(U_{l}-U_{l-1}) \|_{3}.\notag
\end{align}
Hence, we obtain
\begin{align}
	\| \partial_{i}(U_{l+1}-U_{l}) \|_{3}
	\leq C_{7}(l+1)2^{-l(p-1)} 
	\quad \mbox{for} \quad l \geq 1.\label{dai32}
\end{align}
The estimates (\ref{dai18}) and (\ref{dai32}) imply that the $\nabla_x^{\alpha}U_{l}$ for 
$|\alpha| \leq 1$ converge uniformly for  $l \to \infty$
towards functions $\nabla_x^{\alpha} U$, which are continuous in
$x$, $t$, where $U$ is a solution of (\ref{dai1}).

We put
\begin{align}
	C\e_{0}^{p(p-1)} 
	6^{\gamma(p,2)/2}\leq 1.\nonumber
\end{align}
For $0<\e \leq \e_0$, if we assume that
\begin{align}
	C \e^{p(p-1)}
	\left(\frac{4T}{k}\right)^{\gamma(p,2)/2}\leq1,\nonumber
\end{align}
then (\ref{dai3}) holds. Hence, (\ref{lower_lifespan}) in the case of 
$\d\int_{\R^2} g(x) dx=0$ is proved for $0 < \e \leq \e_{0}$.

Next, we shall show the estimate for (\ref{lower_lifespan}) in the case of
\newline $\d\int_{\R^2} g(x) dx \neq 0$.
Let $Y$ be the linear space defined by 
\begin{align}
	Y=&\{u(x,t) \ : \ \nabla_x^{\alpha}u(x,t) \in C(\R^2 \times [0,T]),
	\ \| \nabla_x^{\alpha} u \|_{1} < \infty \quad 
	\mbox{for} \quad |\alpha| \leq 1,\ \nonumber\\
		 &\mbox{supp} \ u
	\subset \{(x,t):|x| \leq t+k \}	
	\}.\nonumber
\end{align}
We can verify that $Y$ is complete with respect to
the norm
\begin{align}
	\| u \|_{Y} =
	 \sum_{|\alpha| \leq 1} \| \nabla_{x}^{\alpha} u \|_1.\nonumber
\end{align}
We define the sequence of functions $\{ u_{l}\}$ by
\[
	u_{1}=u^{0},\quad u_{l}=u^{0}+L(|u_{l-1}|) \ \mbox{for}
	\ l \geq 2.
\]
Since Lemma \ref{lem:decay_est_v} yields
that $\|u^{0}\|_{Y} \leq C'_{0}\e$ where
\[
C'_0:=\wt{C_0}\left|\int_{\R^2}g(x)dx\right|+k^{-1}C_0,
\]
we have that  $u^{0} \in Y$.
We now assume that $\e$ is so small that
\begin{align}
	2^p p C_{1}k^2(C'_{0}\e)^{p-1}D_{1}(T) \leq 1.\label{dai42}
\end{align}
Therefore, as in \cite{J79}, we see that if (\ref{dai42}) holds, then 
there exists a unique local solution of (\ref{IE}). 
By taking $\e_{0}$ small,
the lower bound (\ref{lower_lifespan})
with $\d\int_{\R^2} g(x) dx \neq 0$ follows immediately from (\ref{dai42}).
\hfill$\Box$

\newpage
\par\noindent
{\bf Acknowledgement}
\par
The third author is partially supported
by the Grant-in-Aid for Scientific Research (C) (No.15K04964), 
Japan Society for the Promotion of Science,
and Special Research Expenses in FY2016, General Topics (No.B21),
Future University Hakodate.
All the authors are grateful to the referee for his or her
pointing out many typos in the manuscript.


\end{document}